\documentclass{article}
\usepackage{graphicx}
\usepackage{url}
\usepackage{latexsym}
\usepackage[dvips]{color}
\usepackage[font=small, margin=.75cm]{caption}
\usepackage{amsbsy,amsmath,amsfonts,amsthm,dsfont,latexsym,amssymb,alltt,array}
\usepackage{hyperref}

\newcommand{\pitwo}{\Pi^{\mathrm{P}}_2}

\newcommand{\bbZ}{\mathds{Z}}
\newcommand{\bbR}{\mathds{R}}

\newcommand{\arw}{\rightarrow}
\newcommand{\tri}{\triangle}

\newcommand{\abs}[1]{\lvert #1 \rvert} 
\newcommand{\norm}[1]{\lVert #1 \rVert} 

\theoremstyle{plain}
\newtheorem{theorem}{Theorem}

\theoremstyle{definition}

\theoremstyle{remark}

\theoremstyle{remark}

\theoremstyle{remark}

\makeatletter
\newfont{\footsc}{cmcsc10 at 8truept}
\newfont{\footbf}{cmbx10 at 8truept}
\newfont{\footrm}{cmr10 at 10truept}

\makeatother

\title {\bf Use of MAX-CUT for \\ Ramsey Arrowing of Triangles
}
\author{
Alexander R. Lange \\ Stanis\l aw P. Radziszowski \\
\small Department of Computer Science \\[-0.8ex]
\small Rochester Institute of Technology \\[-0.8ex]
\small Rochester, NY  14623 \\[-0.8ex]
\small \texttt{\{arl9577,spr\}@cs.rit.edu}\\
\\ [-.6em]
and\\
\\ [-.6em]
Xiaodong Xu\\
\small Guangxi Academy of Sciences\\[-0.8ex]
\small Nanning, Guangxi 530007, China\\[-0.8ex]
\small \texttt{xxdmaths@sina.com}\\
}

\date{}

\begin{document}

\maketitle

\vspace*{.2em}

\begin{abstract}
In 1967, Erd\H{o}s and Hajnal asked the question: Does there exist a $K_4$-free
graph that is not the union of two triangle-free graphs?  Finding such a graph
involves solving a special case of the classical Ramsey arrowing
operation. Folkman proved the existence of these graphs in 1970, and they are
now called Folkman graphs.  Erd\H{o}s offered \$100 for deciding if one exists
with less than $10^{10}$ vertices. This problem remained open until 1988 when
Spencer, in a seminal paper using probabilistic techniques, proved the existence
of a Folkman graph of order $3\times 10^9$ (after an erratum), without
explicitly constructing it. In 2008, Dudek and R\"{o}dl developed a strategy to
construct new Folkman graphs by approximating the maximum cut of a related
graph, and used it to improve the upper bound to 941. We improve this bound
first to 860 using their approximation technique and then further to 786 with
the MAX-CUT semidefinite programming relaxation as used in the
Goemans-Williamson algorithm.
\end{abstract}


\bigskip
\section{Introduction}
Given a simple graph $G$, we write $G \arw (a_1,\dots,a_k)^e$ and say that $G$
\emph{arrows} $(a_1,\dots,a_k)^e$ if for every edge $k$-coloring of $G$, a
monochromatic $K_{a_i}$ is forced for some color $i \in
\{1,\dots,k\}$. Likewise, for graphs $F$ and $H$, $G\arw (F,H)^e$ if for every
edge 2-coloring of $G$, a monochromatic $F$ is forced in the first color or a
monochromatic $H$ is forced in the second. Define
$\mathcal{F}_e(a_1,\dots,a_k;p)$ to be the set of all graphs that arrow
$(a_1,\dots,a_k)^e$ and do not contain $K_p$; they are often called Folkman
graphs. The edge Folkman number $F_e(a_1,\dots,a_k;p)$ is the smallest order of
a graph that is a member of $\mathcal{F}_e(a_1,\dots,a_k;p)$. In 1970, Folkman
\cite{Folkman1970} showed that for $k > \max{\{s,t\}}$, $F_e(s,t;k)$ exists. The
related problem of vertex Folkman numbers, where vertices are colored instead of
edges, is more studied \cite{Luczak2001,Nenov2003} than edge Folkman numbers,
but we will not be discussing them. Therefore, we will skip the use of the
superscript $e$ when discussing arrowing, as it is usually used to distinguish
between edge and vertex colorings.

In 1967, Erd\H{o}s and Hajnal \cite{Erdos1967} asked the question: Does there
exist a $K_4$-free graph that is not the union of two triangle-free graphs? This
question is equivalent to asking for the existence of a $K_4$-free graph such
that in any edge 2-coloring, a monochromatic triangle is forced. After Folkman
proved the existence of such a graph, the question then became to find how small
this graph could be, or using the above notation, what is the value of $F_e(3,3;4)$.
Prior to this paper, the best known bounds for this case were $19 \leq
F_e(3,3;4) \leq 941$ \cite{Radziszowski2007, Dudek2008}.

Folkman numbers are related to Ramsey numbers $R(s,t)$, which are defined as the least
positive $n$ such that any 2-coloring of the edges of $K_n$ yields a
monochromatic $K_s$ in the first color or a monochromatic $K_t$ in the
second color. Using the arrowing operator, it is clear that $R(s,t)$ is the smallest
$n$ such that $K_n \arw (s,t)$. The known values and bounds for various types
of Ramsey numbers are collected and regularly updated by the second author
\cite{Radziszowski2011}.

We will be using standard graph theory notation: $V(G)$ and $E(G)$ for the
vertex and edge sets of graph $G$, respectively. A \emph{cut} is a partition of
the vertices of a graph into two sets, $S \subset V(G)$ and
$\overline{S}=V(G)\setminus S$. The \emph{size} of a cut is the number of edges
that join the two sets, that is, $\abs{ \{ \{u,v\} \in E(G)\; | \; u \in S
  \text{ and } v \in \overline{S} \} }$. MAX-CUT is a well-known
\textbf{NP}-hard combinatorial optimization problem which asks for the maximum
size of a cut of a graph.

\newpage

\section{History of $F_e(3,3;4)$} 

\begin{table}[h]
  \centering
  { \renewcommand{\arraystretch}{1.2}
    \begin{tabular}{ l@{$\;\;\;$} | r@{ -- } l| l r }\hline
      Year &  \multicolumn{2}{m{2.7cm}|}{\centering Lower/Upper Bounds} & Who/What & Ref.\\ \hline
      1967 & \multicolumn{2}{c|}{any?$\quad$} & Erd\H{o}s-Hajnal &\cite{Erdos1967} \\
      1970 & \multicolumn{2}{c|}{exist$\quad$} & Folkman &\cite{Folkman1970}\\
      1972 & 10 & & Lin &\cite{Lin1972}\\
      1975 & & $10^{10}$? & Erd\H{o}s offers \$100 for proof &\\
      1986 & & $8 \times 10^{11}$ & Frankl-R\"{o}dl &\cite{Frankl1986}\\
      1988 & & $3\times 10^9$ & Spencer &\cite{Spencer1988} \\
      1999 & $\;\; 16$ & & Piwakowski et al. (implicit) &\cite{Piwakowski1999}\\
      2007 & 19 & & Radziszowski-Xu &\cite{Radziszowski2007}\\
      2008 & & 9697 & Lu &\cite{Lu2008}\\
      2008 & & 941 & Dudek-R\"{o}dl &\cite{Dudek2008} \\
      2012 & & 786 & this work &\\
      2012 & & 100? & Graham offers \$100 for proof &\\ \hline
    \end{tabular}
    }
  \caption{Timeline of progress on $F_e(3,3;4)$.}
  \label{tab:hist}
\end{table}

Table \ref{tab:hist} summarizes the events surrounding $F_e(3,3;4)$, starting
with Erd\H{o}s and Hajnal's \cite{Erdos1967} original question of
existence. After Folkman \cite{Folkman1970} proved the existence, Erd\H{o}s, in
1975, offered \$100 for deciding if $F_e(3,3;4)<10^{10}$. This question remained
open for over 10 years. Frankl and R\"{o}dl \cite{Frankl1986} nearly met
Erd\H{o}s' request in 1986 when they showed that $F_e(3,3;4)$ $< 7.02 \times
10^{11}$. In 1988, Spencer \cite{Spencer1988}, in a seminal paper using
probabilistic techniques, proved the existence of a Folkman graph of order
$3\times 10^9$ (after an erratum by Hovey), without explicitly constructing
it. In 2007, Lu showed that $F_e(3,3;4)\leq 9697$ by constructing a family of
$K_4$-free circulant graphs (which we discuss in Section \ref{sec:lu}) and
showing that some such graphs arrow $(3,3)$ using spectral analysis. Later,
Dudek and R\"{o}dl reduced the upper bound to the best known to date,
$941$. Their method, which we have pursued further with some success, is
discussed in the next section.

The lower bound for $F_e(3,3;4)$ was much less studied than the upper bound. Lin
\cite{Lin1972} obtained a lower bound on $10$ in 1972 without the help of a
computer. All 659 graphs on 15 vertices witnessing $F_e(3,3;5)=15$
\cite{Piwakowski1999} contain $K_4$, thus giving the bound $16 \leq
F_e(3,3;4)$. In 2007, two of the authors of this paper gave a computer-free
proof of $18 \leq F_e(3,3;4)$ and improved the lower bound further to $19$ with
the help of computations \cite{Radziszowski2007}.

The long history of $F_e(3,3;4)$ is not only interesting in itself but also
gives insight into how difficult the problem is. Finding good bounds on the
smallest order of any Folkman graph (with fixed parameters) seems to be
difficult, and some related Ramsey graph coloring problems are \textbf{NP}-hard
or lie even higher in the polynomial hierarchy. For example, Burr
\cite{Burr1976} showed that arrowing $(3,3)$ is $\mathbf{coNP}$-complete, and
Schaefer \cite{Schaefer2001} showed that for general graphs $F$, $G$, and $H$,
$F\arw (G,H)$ is $\mathbf{\pitwo}$-complete.

\section{Arrowing via MAX-CUT}
Building off Spencer's and other methods, Dudek and R\"{o}dl \cite{Dudek2008} in
2008 showed how to construct a graph $H_G$ from a graph $G$, such that the
maximum size of a cut of $H_G$ determines whether or not $G \arw (3,3)$. They
construct the graph $H_G$ as follows. The vertices of $H_G$ are the edges of
$G$, so $\abs{V(H_G)}=\abs{E(G)}$. For $e_1,e_2 \in V(H_G)$, if edges
$\{e_1,e_2,e_3\}$ form a triangle in $G$, then $\{e_1,e_2\}$ is an edge in
$H_G$.

Let $t_\tri(G)$ denote the number of triangles in graph $G$. Clearly,
$\abs{E(H_G)}$$=3t_\tri(G)$. Let $MC(H)$ denote the MAX-CUT value of graph $H$.

\bigskip
\begin{theorem}[Dudek and R\"{o}dl \cite{Dudek2008}] \label{th:mc}
$G \arw (3,3)$ if and only if\\ $MC(H_G) < 2t_\tri(G)$.
\end{theorem}

There is a clear intuition behind Theorem \ref{th:mc} that we will now
describe. Any edge $2$-coloring of $G$ corresponds to a bipartition of the
vertices in $H_G$. If a triangle colored in $G$ is not monochromatic, then its
three edges, which are vertices of $H_G$, will be separated in the bipartition. If
we treat this bipartition as a cut, then the size of the cut will count each
triangle twice for the two edges that cross it. Since there is only one triangle
in a graph that contains two given edges, this effectively counts the number of
non-monochromatic triangles. Therefore, if it is possible to find a cut that has
size equal to $2t_\tri(G)$, then such a cut defines an edge coloring of $G$ that
has no monochromatic triangles. However, if $MC(H_G)<2t_\tri(G)$, then in each
coloring, all three edges of some triangle are in one part and thus, $G \arw
(3,3)$.

A benefit of converting the problem of arrowing $(3,3)$ to MAX-CUT is that the
latter is well-known and has been studied extensively in computer science and
mathematics (see for example \cite{Commander2009}). The decision problem
MAX-CUT$(H,k)$ asks whether or not $MC(H)\geq k$. It is known that MAX-CUT is
\textbf{NP}-hard and this decision problem was one of Karp's 21
\textbf{NP}-complete problems \cite{Karp1972}. In our case, $G\arw (3,3)$ if
and only if MAX-CUT$\left(H_G,2t_\tri(G)\right)$ doesn't hold. Since MAX-CUT is
\textbf{NP}-hard, an attempt is often made to approximate it, such as in the
approaches presented in the next two sections.

\subsection{Minimum Eigenvalue Method}
A method exploiting the minimum eigenvalue was used by Dudek and R\"{o}dl
\cite{Dudek2008} to show that some large graphs are members of
$\mathcal{F}_e(3,3;4)$. The following upper bound \eqref{eq:mineig} on $MC(H_G)$
can be found in \cite{Dudek2008}, where $\lambda_{\text{min}}$ denotes the
minimum eigenvalue of the adjacency matrix of $H_G$.

\begin{equation}\label{eq:mineig}
  MC(H_G) \leq \frac{\abs{E(H_G)}}{2} - \frac{\lambda_{\text{min}}\abs{V(H_G)}}{4}.
\end{equation}

For positive integers $r$ and $n$, if $-1$ is an $r$-th residue modulo $n$, then
let $G(n,r)$ be a circulant graph on $n$ vertices with the vertex set $\bbZ_n$
and the edge set $E(G(n,r)) = \{ \{u,v\} \; | \; u \neq v \text{ and } u-v
\equiv \alpha^r \bmod{n}, \text{ for some } \alpha \in \bbZ_n \}$.

The graph $G_{941}=G(941,5)$ has 707632 triangles. Using the MATLAB
\cite{MATLAB2011} {\tt eigs} function, Dudek and R\"{o}dl \cite{Dudek2008}
computed
\begin{equation*}
  MC(H_{G_{941}})\leq 1397484 < 1415264=2t_\tri(G_{941}).
\end{equation*}
 Thus, by Theorem 1, $G_{941}$ $\arw$ $(3,3)$.
\vspace*{1em}

In an attempt to improve $F_e(3,3;4) \leq 941$, we tried removing vertices
of $G_{941}$ to see if the minimum eigenvalue bound would still show
arrowing. We applied multiple strategies for removing vertices, including removing neighborhoods of vertices, randomly selected vertices, and independent sets of vertices. Most of these strategies were successful, and led to the following theorem:

\bigskip
\begin{theorem}\label{th:gc}
  $F_e(3,3;4) \leq 860$.
\end{theorem}
\vspace*{.5em}

\noindent \textbf{Proof.} For a graph $G$ with vertices $\bbZ_n$, define
$C=C(d,k) = \{ v \in V(G)\;|\; v=id \bmod{n}, \text{ for } 0 \leq i < k\}$. Let
$G=G_{941}$, $d=2$, $k=81$, and $G_C$ be the graph induced on $V(G) \setminus
C(d,k)$. Then $G_C$ has 860 vertices, 73981 edges and 542514 triangles. Using
the MATLAB {\tt eigs} function, we obtain $\lambda_{\text{min}} \approx
-14.663012$. Setting $\lambda_{\text{min}} > -14.664$ in \eqref{eq:mineig} gives
\begin{equation}
MC(H_{G_C}) < 1084985 < 1085028=2t_\tri(G_C).
\end{equation}

\noindent Therefore, $G_C \arw (3,3)$. $\Box$
\vspace*{1em}

None of the methods used allowed for $82$ or more vertices to be removed without
the upper bound on $MC$ becoming larger than $2t_\tri$.

\subsection{Goemans-Williamson Method}
The Goemans-Williamson MAX-CUT approximation algorithm \cite{Goemans1995} is a
well-known, polynomial-time algorithm that relaxes the problem to a semi-definite 
program (SDP). It involves the first use of SDP in combinatorial approximation
and has since inspired a variety of other successful algorithms (see for example
\cite{Karloff1997,Frieze1997}). This randomized algorithm returns a cut with
expected size at least 0.87856 of the optimal value. However, in our case, all
that is needed is a feasible solution to the SDP, as it gives an upper bound on
$MC(H)$. A brief description of the Goemans-Williamson relaxation follows.

The first step in relaxing MAX-CUT is to represent the problem as a quadratic
integer program. Given a graph $H$ with $V(H)=\{1,\dots,n\}$ and nonnegative
weights $w_{i,j}$ for each pair of vertices $\{i,j\}$, we can write $MC(H)$ as
the following objective function:
  \begin{align}
    \text{Maximize}\quad &\frac{1}{2}\sum_{i<j}w_{i,j}(1 -
    y_iy_j) \label{eq:quadratic}\\ \text{subject to:}\quad &y_i \in \{-1,1\}
    \quad \text{for all } i \in V(H). \notag
  \end{align}
 
Define one part of the cut as $S=\{i\; |\; y_i = 1\}$. Since in our case all
graphs are weightless, we will use
\begin{equation*}
w_{i,j} =
\begin{cases}
  1 & \text{if } \{i,j\} \in E(H),\\
  0 & \text{otherwise}.
\end{cases}
\end{equation*}

Next, the integer program \eqref{eq:quadratic} is relaxed by extending the
problem to higher dimensions. Each $y_i \in \{-1,1\}$ is now replaced with a
vector on the unit sphere $\mathbf{v}_i \in \bbR^n$, as follows:

\begin{align}
  \text{Maximize}\quad &\frac{1}{2}\sum_{i<j}w_{i,j}(1 - \mathbf{v}_i\cdot
  \mathbf{v}_j) \label{eq:relax}\\ \text{subject to:}\quad &\norm{\mathbf{v}_i}
  = 1 \quad \text{for all } i \in V(H). \notag
\end{align}

If we define a matrix $Y$ with the entries
$y_{i,j}=\mathbf{v_i}\cdot\mathbf{v_j}$, that is, the Gram matrix of ${
  \mathbf{v}_1,\dots,\mathbf{v}_n}$, then $y_{i,i}=1$ and $Y$ is positive
semidefinite. Therefore, \eqref{eq:relax} is a semidefinite program.

\vspace*{1em}

\subsection{Some Cases of Arrowing} \label{sec:lu}

Using the Goemans-Williamson approach, we tested a wide variety of graphs for
arrowing by finding upper bounds on MAX-CUT. These graphs included the $G(n,r)$
graphs tested by Dudek and R\"{o}dl, similar circulant graphs based on the Galois fields $GF(p^k)$, and random graphs. Various modifications of these graphs
were also considered, including the removal and/or addition of vertices and/or
edges, as well as copying or joining multiple candidate graphs together in
various ways. We tested the graph $G_C$ of Theorem \ref{th:gc} and obtained the
upper bound $MC(H_{G_C}) \leq 1077834$, a significant improvement over the bound
$1084985$ obtained from the minimum eigenvalue method. This provides further
evidence that $G_C \arw (3,3)$, and is an example of when
\eqref{eq:relax} yields a much better upper bound.

Multiple SDP solvers that were designed \cite{Burer2003,Helmberg2000} to handle
large-scale SDP and MAX-CUT problems were used for the tests. Specifically, we
made use of a version of {\tt SDPLR} by Samuel Burer \cite{Burer2003}, a solver
that uses low-rank factorization. The version {\tt SDPLR-MC} includes
specialized code for the MAX-CUT SDP relaxation. {\tt SBmethod} by Christoph
Helmberg \cite{Helmberg2000} implements a spectral bundle method and was also
applied successfully in our experiments. In all cases where more than one solver
was used, the same results were obtained.

The type of graph that led to the best results was described by Lu
\cite{Lu2008}.  For positive integers $n$ and $s$, $s<n$, $s$ relatively prime
to $n$, define set $S = \{ s^i \bmod{n} \; | \; i=0,1,\dots,m-1\}$, where $m$ is
the smallest positive integer such that $s^m \equiv 1 \bmod{n}$. If $-1 \bmod{n}
\in S$, then let $L(n,s)$ be a circulant graph on $n$ vertices with
$V(L(n,s))=\bbZ_n$. For vertices $u$ and $v$, $\{u,v\}$ is an edge of $L(n,s)$
if and only if $u-v \in S$. Note that the condition that $-1 \bmod{n} \in S$
implies that if $u-v \in S$ then $v-u \in S$.

In Table 1 of \cite{Lu2008}, a set of potential members of
$\mathcal{F}_e(3,3;4)$ of the form $L(n,s)$ were listed, and the graph
$L(9697,4)$ was shown to arrow $(3,3)$.  Lu gave credit to Exoo for showing
that $L(17,2)$, $L(61,8)$, $L(79,12)$, $L(421,7)$, and $L(631,24)$ do not arrow
$(3,3)$.
\vspace*{1em}

 We tested all graphs from Table 1 of \cite{Lu2008} of order less than 941 with
the MAX-CUT method, using both the minimum eigenvalue and SDP upper
bounds. Table \ref{tab:results} lists the results. Note that although none of
the computed upper bounds of the $L(n,s)$ graphs imply arrowing $(3,3)$, all SDP
bounds match those of the minimum eigenvalue bound. This is distinct from other
families of graphs, including those in \cite{Dudek2008}, as the SDP bound is
usually tighter. Thus, these graphs were given further consideration.

\begin{table}[h]
  \centering
  { \renewcommand{\arraystretch}{1.2}
  \begin{tabular}{ | c | r | r | r | r | r | r | } \hline
    $G$ & \multicolumn{1}{c|}{$2t_\tri(G)$} &
    \multicolumn{1}{c|}{$\lambda_{\text{min}}$} &
    \multicolumn{1}{c|}{SDP}\\ \hline\hline $L(127,5)$ & 19558 & 20181 & 20181
    \\ $L(457,6)$ & 347320 & 358204 & 358204 \\ $L(761,3)$ & 694032 & 731858 &
    731858 \\ $L(785,53)$ & 857220 & 857220 & 857220 \\ \hline $G_{786}$ &
    857762 & 857843 & 857753 \\ \hline
  \end{tabular}
  }
  \caption{Potential $\mathcal{F}_e(3,3;4)$ graphs $G$ and upper bounds on
    $MC(H_G)$, where ``$\lambda_{\text{min}}$'' is the bound \eqref{eq:mineig}
    and ``SDP'' is the solution of \eqref{eq:relax} from {\tt SDPLR-MC} and {\tt
      SBmethod}. $G_{786}$ is the graph of Theorem \ref{th:786}.}
  \label{tab:results}
\end{table}

$L(127,5)$ was given particular attention, as it is the same graph as $G_{127}$,
where $V(G_{127})=\bbZ_{127}$ and $E(G_{127})=\{ \{x,y\} \; | \; x-y \equiv
\alpha^3 \bmod{127} \}$ (that is, the graph $G(127,3)$ as defined in the
previous section). It has been conjectured by Exoo that $G_{127} \arw
(3,3)$. He also suggested that subgraphs induced on less than 100 vertices of
$G_{127}$ may as well. For more information on $G_{127}$ see
\cite{Radziszowski2007}.

Numerous attempts were made at modifying these graphs in hopes that one of the
MAX-CUT methods would be able to prove arrowing. Indeed, we were able to do so
with $L(785,53)$. Notice that all of the upper bounds for $MC(H_{L(785,53)})$
are $857220$, the same as $2t_\tri\left(L(785,53)\right)$. Our goal was then to
slightly modify $L(785,53)$ so that this value becomes smaller. Let $G_{786}$
denote the graph $L(785,53)$ with one additional vertex connected to the
following 60 vertices:

\vspace*{1em}

\begin{minipage}{\textwidth} \centering
{
\small
\begin{verbatim}
 {  0,   1,   3,   4,   6,   7,   9,  10,  12,  13,  15,  16,
   18,  19,  21,  22,  24,  25,  27,  28,  30,  31,  33,  34,
   36,  37,  39,  40,  42,  43,  45,  46,  48,  49,  51,  52,
   54,  55,  57,  58,  60,  61,  63,  66,  69, 201, 204, 207,
  210, 213, 216, 219, 222, 225, 416, 419, 422, 630, 642, 645  }
\end{verbatim}
}
\end{minipage}

\vspace*{1em}


 $G_{786}$ is still $K_4$-free, has 61290 edges, and has 428881 triangles. The
upper bound computed from the SDP solvers for $MC(H_{G_{786}})$ is 857753. We
did not find a nice description for the vectors of this solution.  Software
implementing {\tt SpeeDP} by Grippo et al. \cite{Grippo2010b}, an algorithm
designed to solve large MAX-CUT SDP relaxations, was used by Rinaldi (one of the
authors of \cite{Grippo2010b}) to analyze this graph. He was able to obtain the
bounds $857742 \leq MC(H_{G_{786}}) \leq 857750$, which agrees with, and
improves over our upper bound computation. Since $2t_\tri(G_{786}) = 857762$, we
have both from our tests and his {\tt SpeeDP} test that $G_{786} \arw (3,3)$,
and the following main result.

\bigskip
\begin{theorem}\label{th:786}
$F_e(3,3;4) \leq 786.$
\end{theorem}

We note that finding a lower bound on MAX-CUT, such as the $857742 \leq MC(H_{G_{786}})$ bound from {\tt SpeeDP}, follows from finding
an actual cut of a certain size. This method may be useful, as
finding a cut of size $2t_\tri(G)$ shows that $G \not\arw (3,3)$.

\section{Tasks to Complete}
Improving the upper bound on $F_e(3,3;4)$ $\leq 786$ is the main challenge. The
question of whether $G_{127} \arw (3,3)$ is still open, and any method that
could solve it would be of much interest.

During the 2012 SIAM Conference on Discrete Mathematics in Halifax, Nova Scotia,
Ronald Graham announced a \$100 award for determining if $F_e(3,3;4) < 100$.

Another open question is the lower bound on $F_e(3,3;4)$, as it is quite
puzzling that only 19 is the best known. Even an improvement to $20 \leq
F_e(3,3;4)$ would be good progress.

\section{Acknowledgments}
The third author is supported by the Guangxi Natural Science Foundation
(2011GXNSFA018142). We would like to thank Giovanni Rinaldi and Luigi Grippo
for their enthusiastic aid in the computation of MAX-CUT bounds with their {\tt
  SpeeDP} algorithm \cite{Grippo2010b}. We would also like to thank the referee for the 
helpful comments.
\bibliographystyle{plain}
\bibliography{fe334.bib}

\begin{thebibliography}{10}

\bibitem{Burer2003}
Samuel Burer and Renato~D.C. Monteiro.
\newblock {A nonlinear programming algorithm for solving semidefinite programs
  via low-rank factorization}.
\newblock {\em Mathematical Programming (Series B)}, 95(2):329--357, February
  2003.
\newblock Software available at {\tt http://dollar.biz.uiowa.edu/\~{}sburer}.

\bibitem{Burr1976}
Stefan~A. Burr.
\newblock 1976. {R}esult mentioned in book by {M}. {G}arey and {D}. {J}ohnson.
\newblock {\em Computers and Intractability: A Guide to the Theory of
  NP-Completeness}, 1979.
\newblock W. H. Freeman and Company.

\bibitem{Commander2009}
Clayton~W. Commander.
\newblock {Maximum Cut Problem, MAX-CUT}.
\newblock In Christodoulos Floudas and Panos Pardalos, editors, {\em
  Encyclopedia of Optimization}, pages 1991--1999. Springer, second edition,
  2009.

\bibitem{Dudek2008}
Andrzej Dudek and Vojtech R\"odl.
\newblock {On the Folkman Number $f(2,3,4)$}.
\newblock {\em Experimental Mathematics}, 17(1):63--67, 2008.

\bibitem{Erdos1967}
Paul Erd\H{o}s and Andr\'{a}s Hajnal.
\newblock Research problem 2--5.
\newblock {\em Journal of Combinatorial Theory}, 2:104, 1967.

\bibitem{Folkman1970}
Jon Folkman.
\newblock Graphs with monochromatic complete subgraphs in every edge coloring.
\newblock {\em SIAM Journal on Applied Mathematics}, 18(1):19--24, January
  1970.

\bibitem{Frankl1986}
Peter Frankl and Vojtech R\"{o}dl.
\newblock Large triangle-free subgraphs in graphs without ${K}_4$.
\newblock {\em Graphs and Combinatorics}, 2:135--144, 1986.

\bibitem{Frieze1997}
Alan Frieze and Mark Jerrum.
\newblock {Improved Approximation Algorithms for MAX $k$-CUT and MAX
  BISECTION}.
\newblock {\em Algorithmica}, 18(1):67--81, 1997.

\bibitem{Goemans1995}
Michael Goemans and David Williamson.
\newblock {Improved Approximation Algorithms for Maximum Cut and Satisfiability
  Problems Using Semidefinite Programming}.
\newblock {\em Journal of the ACM}, 42(6):1115--1145, 1995.

\bibitem{Grippo2010b}
Luigi Grippo, Laura Palagi, Mauro Piacentini, Veronica Piccialli, and Giovanni
  Rinaldi.
\newblock {SpeeDP}: An algorithm to compute {SDP} bounds for very large
  {Max-Cut} instances.
\newblock {\em Mathematical Programming}, 2012.
\newblock doi:10.1007/s10107-012-0593-0.

\bibitem{Helmberg2000}
Christoph Helmberg and Franz Rendl.
\newblock {A Spectral Bundle Method for Semidefinite Programming}.
\newblock {\em SIAM Journal of Optimization}, 10:673--696, 2000.
\newblock Software available at {\tt
  http://www-user.tu-chemnitz.de/\~{}helmberg}.

\bibitem{Karloff1997}
Howard Karloff and Uri Zwick.
\newblock {A $7/8$ Approximation Algorithm for MAX 3SAT?}
\newblock In {\em 38th Annual IEEE Symposium on Foundations of Computer
  Science}, pages 406--415, 1997.

\bibitem{Karp1972}
Richard~M. Karp.
\newblock {Reducibility Among Combinatorial Problems}.
\newblock In R.~E. Miller and J.~W. Thatcher, editors, {\em Complexity of
  Computer Computations}, pages 85--103. Plenum, New York, 1972.

\bibitem{Lin1972}
Shen Lin.
\newblock On {R}amsey numbers and ${K}_r$-coloring of graphs.
\newblock {\em Journal of Combinatorial Theory, Series B}, 12:82--92, 1972.

\bibitem{Lu2008}
Linyuan Lu.
\newblock Explicit {C}onstruction of {S}mall {F}olkman {G}raphs.
\newblock {\em SIAM Journal on Discrete Mathematics}, 21(4):1053--1060, January
  2008.

\bibitem{Luczak2001}
Tomasz {\L}uczak, Andrzej Ruci\'{n}ski, and Sebastian Urba\'{n}ski.
\newblock {On minimal Folkman graphs}.
\newblock {\em Discrete Mathematics}, 236:245--262, 2001.

\bibitem{MATLAB2011}
MATLAB.
\newblock {\em \textnormal{{V}ersion 7.12.0 (R2011a)}}.
\newblock The MathWorks Inc., Natick, Massachusetts, 2011.
\newblock \\{\tt http://www.mathworks.com/products/matlab}.

\bibitem{Nenov2003}
Nedyalko Nenov.
\newblock {On the triangle vertex Folkman numbers}.
\newblock {\em Discrete Mathematics}, 271:327--334, September 2003.

\bibitem{Piwakowski1999}
Konrad Piwakowski, Stanis{\l}aw~P. Radziszowski, and Sebastian Urba\'{n}ski.
\newblock {Computation of the Folkman Number $F_e (3,3;5)$}.
\newblock {\em Journal of Graph Theory}, 32:41--49, 1999.

\bibitem{Radziszowski2011}
Stanis{\l}aw~P. Radziszowski.
\newblock {Small Ramsey Numbers}.
\newblock {\em Electronic Journal of Combinatorics}, August 2011.
\newblock Dynamic Survey 1, Revision
  \#13.$\;\;$\url{http://www.combinatorics.org}.

\bibitem{Radziszowski2007}
Stanis{\l}aw~P. Radziszowski and Xiaodong Xu.
\newblock {On the Most Wanted Folkman Graph}.
\newblock {\em Geombinatorics}, 16(4):367--381, 2007.

\bibitem{Schaefer2001}
Marcus Schaefer.
\newblock {Graph Ramsey Theory and the Polynomial Hierarchy}.
\newblock {\em Journal of Computer and System Sciences}, 62:290--322, 2001.

\bibitem{Spencer1988}
Joel Spencer.
\newblock Three hundred million points suffice.
\newblock {\em Journal of Combinatorial Theory, Series A}, 49(2):210--217,
  1988.
\newblock Also see erratum by M. Hovey in Vol. 50, p. 323.

\end{thebibliography}

\end{document}